\documentclass[11pt]{amsart}
\usepackage{amssymb,amsfonts}
\usepackage{euscript,mathrsfs}
\usepackage{latexsym}
\usepackage{amscd}
\usepackage{amsmath,bbm}
\usepackage{a4wide}
\usepackage{epsf}
\usepackage{color}
\usepackage{graphicx}
\usepackage[utf8]{inputenc}
\usepackage{amsmath,amsthm,amsopn,fixltx2e,mparhack,graphicx,microtype}
\usepackage{url}

\newcommand\mm{{\mathfrak m}}

\renewcommand{\>}{\rangle}
\def\ol#1{{\overline {#1}}}
\newcommand\oJ{{\hspace{.45ex}\overline{\hspace{-.45ex}J}}}
\DeclareMathOperator{\Hull}{Hull} % The radical
\newcommand{\comment}[1]{ }
\newcommand{\kk}{\mathbbm{k}}

\newcommand{\QQ}{\mathbb{Q}}
\newcommand{\NN}{\mathbb{N}}

\newcommand{\ZZ}{\mathbb{Z}}

\renewcommand{\>}{\rangle}

\renewcommand{\subset}{\subseteq}  % Subsetsymbols

\newcommand{\defas}{\mathrel{\mathop{:}}=}   % Definition

\newcommand{\Mac}{Macaulay~2}
\newcommand{\emb}{\mathrm{emb}}
\newcommand{\cyclo}{\textsl{Cyclotomic}}
\newcommand{\Binom}{\textsl{Binomials}}
\theoremstyle{plain}% default
\newtheorem{thm}{Theorem}

\newtheorem{prop}[thm]{Proposition}

\newtheorem{algorithm}{Algorithm}

\theoremstyle{definition}
\newtheorem{defn}[thm]{Definition}

\newtheorem{exmp}[thm]{Example}

\theoremstyle{remark}

\DeclareMathAlphabet\mathbit
  \encodingdefault\rmdefault\bfdefault\itdefault

\begin{document}

\title[Binomials.m2]{Decompositions of Binomial Ideals in Macaulay 2} 
\author{Thomas Kahle}
\email{www.thomas-kahle.de}

\begin{abstract}
  The package \Binom\ contains implementations of specialized algorithms for binomial ideals,
  including primary decomposition into binomial ideals.  The current implementation works in
  characteristic zero.  Primary decomposition is restricted to binomial ideals with trivial
  coefficients to avoid computations over the algebraic numbers.  The basic ideas of the algorithms
  go back to Eisenbud and Sturmfels' seminal paper on the subject.  Two recent improvements of the
  algorithms are discussed and examples are presented.
\end{abstract}

\maketitle

\section{Binomial ideals}
Let $S = \kk[x_{1},\dots,x_{n}]$ denote the standard polynomial ring over a field $\kk$.  A binomial
ideal $I\subset S$ is an ideal generated by binomials $x^{u}-\lambda x^{v}$, where $u,v\in \NN^{n}$
are exponent vectors and $\lambda\in \kk$ is a coefficient.  Monomials are also considered
binomials.  Assumptions on $\kk$ will be forced upon us when computing primary decompositions.  The
ideal $\<x^{3}-1\>$ has no primary decomposition into binomial ideals when $\kk$ does not contain a
third root of unity.  Interest in binomial ideals is due to the frequency with which they arise in
applications.  To name one, in algebraic statistics one is interested in primary decompositions of
\emph{conditional independence ideals} whose components describe various combinatorial ways in which
a set of conditional independence statements can be
realized~\cite{fink09,juergen09:_binom_edge_ideal_ci}.  Because the minimal primes of binomial
ideals are toric ideals~\cite{eisenbud96:_binom_ideal}, binomial conditional independence models are
unions of exponential families.  In particular they are unirational.  Knowledge of a primary
decomposition also gives a piecewise parameterization of such models.

The new \Mac~\cite{M2} package \Binom\ offers specialized implementations of primary decomposition,
radical computations and minimal and associated primes.  The starting point for this implementation
was Section 9 in Eisenbud and Sturmfels' foundational paper~\cite{eisenbud96:_binom_ideal}, but
various improvements have been discovered and implemented.  \Binom\ is the fastest and often only
way to compute large primary decompositions of binomial ideals.
\begin{exmp}$ $\\[-.4cm]
\begin{verbatim}
i1 : needsPackage "Binomials"
i2 : R = QQ[x,y]
i3 : I = ideal (x^2-x*y, x*y-y^2)

i4 : binomialPrimaryDecomposition I
[...] 
                               2
o4 = {ideal(x - y), ideal (x, y )}
\end{verbatim}
\end{exmp}
A binomial primary decomposition starts with a cellular decomposition.  Recall that a binomial ideal
$I\subset S$ is \emph{cellular} if in $S/I$ every monomial is either regular (i.e.~a nonzerodivisor)
or nilpotent.  The implemented algorithm to compute a cellular decomposition is discussed
in~\cite{kahle09:_decom_binom_ideal,ojeda00:_cellul_binom_ideal}.  Since cellular decomposition is
independent of $\kk$, it can serve as a first approximation of primary decomposition over any field.

In this paper we focus on decomposing a cellular binomial ideal further.  To this end, assume that
$I$ is $J$-cellular for some $J\subset [n]$, that is, the variables with indices in $J$ are regular,
while the variables with indices in $\ol J \defas [n] \setminus J$ are nilpotent.  
% Note that every cellular ideal is $J$-cellular for some $J$.

\section{Computing associated primes}
If $\kk$ is algebraically closed, then the associated primes of a binomial ideal are guaranteed to
be binomial.  Since computer algebra system usually don't implement algebraically closed fields, the
input binomial ideals are restricted to be generated by unital binomials $x^{u}-x^{v}$.  In this
case the binomial primary decomposition together with the associated primes exist over a cyclotomic
extension of~$\QQ$~\cite{kahle09:_decom_binom_ideal}.  If necessary, \Binom\ will construct this
extension and return its result over a different ring.
\begin{exmp}$ $\\[-.4cm]
\begin{verbatim}
i1 : R = QQ[x]
i2 : I = ideal(x^3-1)
i3 : BPD I
[...]
o3 = {ideal(x - 1), ideal(x - ww ), ideal(x + ww  + 1)}
                                3               3
\end{verbatim}
\end{exmp}
In the following discussion we will assume $\kk$ to be algebraically closed and of characteristic
zero.  Let $I$ be $J$-cellular for $J\subset [n]$, and denote $\mm_{J} \defas \<x_{i} : i \notin
J\>$.  The associated primes of $I$ are of the form $I_{\rho, J} + \mm_{J}$, where $I_{\rho,J}
\defas \<x^{u}-\rho(u-v)x^{v} : u-v \in L\>$ is a lattice ideal in the $J$-variables, and $\rho : L
\to \kk^{*}$ is a group homomorphism from a sublattice $L\subset \ZZ^{J}$.  The pair $(\rho, L)$ is
called a \emph{partial character}
in~\cite{dickenstein08:_combin_binom_primar_decom,eisenbud96:_binom_ideal}.  Here we will simply
speak of a character.  An \emph{extension} of $\rho$ is a character $\tau : L' \to \kk$ such that
$L\subset L'$ and $\rho$ and $\tau$ agree on $L$.  A character is \emph{saturated} if its domain is
\emph{saturated lattice}, that is $L$ is not contained with finite index in any other sublattice
of~$\ZZ^{J}$.  An extension to a saturated character is a \emph{saturation}.  Denote $\kk[J] \defas
\kk[x_{i} : i \in J]$, and $\kk[\oJ]\defas\kk[x_{i}:i\notin J]$.  Associated primes of cellular
binomial ideals come in groups.  The following theorem states that they are to be found among the
associated primes of lattice ideals in~$\kk[J]$.
\begin{thm}[{\cite[Theorem~8.1]{eisenbud96:_binom_ideal}}]
  \label{sec:associated-primes-cellular}
  Let $I \subset S$ be a $J$-cellular binomial ideal.  Let $I_{\sigma,J} + \mm_{J}$ be an
  associated prime of~$I$, then there exists a monomial $m \in \kk[\oJ]$, and a character $\tau$ on
  $\ZZ^{J}$ whose saturation is $\sigma$, such that
  \begin{equation*}
    (I: m) \cap \kk[J] = I_{\tau}.
  \end{equation*}
\end{thm}
It can be seen that the converse also holds.  Every associated prime of any occurring lattice ideal
is associated to~$I$~\cite{kahle11mesoprimary}.  Theorem~\ref{sec:associated-primes-cellular} shows
that a sub-problem in the computation of associated primes is to determine the set of lattice ideals
of the form $(I:m) \cap \kk[J]$.
\begin{defn}[\cite{kahle11mesoprimary}]
  A lattice $L \subset \ZZ^{J}$ is \emph{potentially associated} to $I$ if there exists a
  \emph{witness monomial} $m\in \kk[\oJ]$ such that $(I:m) \cap \kk[J] = I_{\rho, J}$ for some
  character $\rho : L \to \kk^{*}$.
\end{defn}
Properly defining the set of associated lattices, which is contained in the set of potentially
associated lattices, requires care and is one of the topics of~\cite{kahle11mesoprimary}.  For
computational purposes the few lattices that are potentially associated, but not associated, play a
minor role.  They will eventually yield redundant primary decompositions, a problem that has to be
handled in any case, since cellular decomposition introduces redundancy on a large scale.  The
lattice ideals $(I:m) \cap \kk[J]$ are partially ordered by inclusion, and so are their lattices.
\begin{defn}
  A potentially associated lattice is called \emph{embedded}, if it properly contains the lattice of
  $I \cap \kk[J]$.
\end{defn}

% \section{Two improvements to binomial primary decomposition algorithms}
% \label{sec:two-impr-binom}
% In this section we will present two improvements to the binomial primary decomposition algorithm in
% \cite{eisenbud96:_binom_ideal,kahle09:_decom_binom_ideal}.  Both of these improvements have been
% implemented on the \Mac-Workshop in G\"{o}ttingen in March 2011.

% \subsection{Searching associated primes}
% \label{sec:search-assoc-prim}
% The first improvement concerns the search for associated primes.
% Theorem~\ref{sec:associated-primes-cellular} tells us how to locate associated primes.  We need to
% examine all nilpotent monomials $m$ of $S/I$ and compute the colon ideals $(I:m)$.  Each of them

A first algorithm to find potentially associated lattices would examine all ideals $(I:m)$ where $m$
is a nonzero monomial in $\kk[\oJ]/ (I \cap \kk[\oJ])$.  By cellularity of $I$ there are only
finitely many such monomials and this search will terminate.  The associated primes algorithm in
\Binom\ instead uses a random search.  The set of monomials to be examined can be very large
compared to relatively few potentially associated lattices.  The design goal in
Algorithm~\ref{alg:random_search} is to compute as few colon ideals $(I:m)$ as possible.  If a
monomial $m$ divides a monomial $n$, then $(I:m) \subset (I:n)$ and containment also holds for the
potentially associated lattices.  Due to this fact we can exclude large posets of monomials if we
find two monomials with the same potentially associated lattice.
\begin{algorithm}\label{alg:random_search}$ $\\
  Input: A $J$-cellular binomial ideal $I$.\\
  Output: The potentially associated lattices of $I$\\
  \begin{enumerate}
  \item Compute the lattice ideal $I\cap \kk[J]$.
  \item Initialize a list of known potentially associated lattices and witnesses containing only the
    pair $(I \cap \kk[J], 1)$.
  \item Initialize a todo-list with all monomials in a $\kk$-basis of $\kk[\oJ]/(I\cap \kk[\oJ])$.
  \item Iterate the following until the todo-list is empty
    \begin{itemize}
    \item Choose and remove a random monomial $m$ from the todo-list.  Compute the lattice ideal
      $(I:m) \cap \kk[J]$ and check if its lattice is already on the list of potentially associated
      lattices.
      \begin{itemize}
      \item   If yes, then add $m$ as a new witness for that lattice, remove from the todo-list
        every monomial between existing witnesses and $m$.
      \item  If no, then add $(I\cap\kk[J], m)$ to the
        list of potentially associated lattices.
      \end{itemize}
    \end{itemize}
  \end{enumerate}
\end{algorithm}
To save space and time, the implementation in \Binom\ does not save all the witness monomials.  If
$m,n$ are both witnesses for the same potentially associated lattice and $m \vert n$, then only $m$
needs to saved.  

Given the set of potentially associated lattices, determining the associated primes is easy.  It
consists of saturating characters and will not be discussed here.  The necessary cyclotomic
extensions are handled in a separate package \cyclo, published together with \Binom.

\section{Computing minimal primary components}
% Computing primary decompositions of binomial ideals relies crucially on the fact that a primary
% binomial ideal contains all binomials of its minimal prime.  The first step in the decomposition
% of a cellular binomial ideal is to compute associated primes using their characterization in
% Theorem~\ref{sec:associated-primes-cellular}.
Let $I\subset S$ be a $J$-cellular binomial ideal, and $P = I_{\rho, J} + \mm_{J}$ one of its
associated primes.  Eisenbud and Sturmfels show that any primary component of $I$ over $P$
contains~$I_{\rho, J}$.  In fact, $I + I_{\rho, J}$ has $P$ as its unique minimal prime and a
primary component is computed by removing all embedded primary components from $I + I_{\rho, J}$.
This is the content of \cite[Theorem~7.1]{eisenbud96:_binom_ideal}.  Let $\Hull(I)$ denote the
intersection of the minimal primary components of a binomial ideal~$I$.  If $I$ is cellular, then
$\Hull(I)$ is binomial.  Computing $\Hull$ of a binomial ideal is a cumbersome procedure.  One way,
described in \cite{eisenbud96:_binom_ideal}, is to successively identify binomials $b$ such that
$(I:b)$ is a binomial ideal strictly containing~$I$.  This approach is slow.  Here we will use a
similar strategy like in Algorithm~\ref{alg:random_search}.
% The search for these binomials and computation of the colon ideal was a serious
% bottleneck in the first implementation of binomial primary decomposition.  Using combinatorial
% information this computation need not be carried out anymore.  For $I$ a $J$-cellular binomial
% ideal, 
Denote $M_{\emb}(I)$ the monomial ideal generated by all witnesses of embedded lattices of~$I$.
Then \cite[Theorem 3.2]{dickenstein08:_combin_binom_primar_decom} implies the following
simplification.
\begin{prop} \label{prop:hull} If $I$ is $J$-cellular and has exactly one minimal prime, then
  \[ \Hull (I) = I + M_{\emb}(I) \] In particular $\Hull(I)$ is binomial.
\end{prop}
To compute the minimal primary component of $I$ over $P = I_{\rho, J} + \mm_{J}$ one computes $\Hull
(I + I_{\rho, J})$~\cite{eisenbud96:_binom_ideal,ojeda00:_cellul_binom_ideal}.  The monomial ideal
$M_{\emb}(I+I_{\rho, J})$ is determined essentially by Algorithm~\ref{alg:random_search}.  It is in
fact simpler, since only minimal generators of $M_{\emb}(I)$ need to be computed.  In most cases
only a small fraction of the standard monomials needs to be examined.
\begin{exmp} This example demonstrates how a component of high multiplicity leads to many monomials
  to be examined.
\begin{verbatim}
i1 : R = QQ[a,b]
i2 : I = ideal (a^10000 * (b-1))
i3 : BPD I
                           10000
o3 = {ideal(b - 1), ideal(a     )}
\end{verbatim}
  In this case the poset of nilpotent monomials is totally ordered and a typical run of
  Algorithm~\ref{alg:random_search} would only compute $\lceil\log_{2}(10000)\rceil = 14$ lattice
  ideals.  Since the structure of the poset of embedded associated lattices can be complicated it is
  not known if there are better search algorithms than random search.
\end{exmp}

\section*{Acknowledgment}
\label{sec:acknowledgement}
The author would like to thank the organizers and participants of \Mac\ workshop March 2011 in
Göttingen on which the improved algorithms were first implemented.

\bibliographystyle{amsplain}
\bibliography{/home/tom/bibtex/math.bib,/home/tom/bibtex/statistics.bib}

\end{document}